# Inequality for Variance of Weighted Sum of Correlated Random Variables and WLLN

Jingwei Liu*

(School of Mathematics and System Sciences, Beihang University, Beijing, 100191, P.R China)

**Abstract**: The upper bound inequality for variance of weighted sum of correlated random variables is derived according to Cauchy-Schwarz 's inequality, while the weights are non-negative with sum of 1. We also give a novel proof with positive semidefinite matrix method. And the variance inequality of sum of correlated random variable with general weights is also obtained. Then, the variance inequalities are applied to the Chebyshev's inequality and sufficient condition of weak law of large numbers (WLLN) for sum of correlated random variables .

## 1. Introduction

Suppose $X_1,\cdots,X_n$ are any random variables (discrete, continuous or else) with finite expectation and variance, for any weights $\alpha_1,\cdots,\alpha_n$ ,which are real numbers. Let $E(X)$ , $Var(X)$ , $Cov(X,Y)$ denote the expectation, variance and covariance respectively. Let

$$\rho=\frac{Cov(X,Y)}{\sqrt{Var(X)}\sqrt{Var(Y)}}$$

denote the correlation coefficient between $X$ and $Y$ , and $|\rho|\le 1$. And, the well-known Cauchy -Schwarz inequality is

$$|Cov(X,Y)|\le\sqrt{Var(X)}\sqrt{Var(Y)}\,. \tag{1}$$

For first and second moments of $\xi=\sum_{i=1}^{n}\alpha_i X_i$ , the expectation is

$$E(\sum_{i=1}^{n}\alpha_i X_i)=\sum_{i=1}^{n}\alpha_i E(X_i)\,. \tag{2}$$

* Corresponding author: jwliu@buaa.edu.cn; jwliu@pku.org.cn; liujingwei03@tsinghua.org.cn

And , the variance of $\xi = \sum_{i=1}^{n} \alpha_i X_i$ is

$$Var(\sum_{i=1}^{n} \alpha_i X_i) = \sum_{i=1}^{n} \alpha_i^2 Var(X_i) + \sum_{1 \le i < j \le n} 2\alpha_i \alpha_j Cov(X_i, X_j) . \quad (3)$$

If the variables are uncorrelated ( a strong condition is independent),

$$Var(\sum_{i=1}^{n} \alpha_i X_i) = \sum_{i=1}^{n} \alpha_i^2 Var(X_i) \quad (4)$$

All of the above definitions and formulae are referred from [1-6]. Motivated by the relationship of $E(\sum_{i=1}^{n} \alpha_i X_i) = \sum_{i=1}^{n} \alpha_i E(X_i)$ , what is the relationship among $Var(\sum_{i=1}^{n} \alpha_i X_i)$ , $\sum_{i=1}^{n} \alpha_i^2 Var(X_i)$ and $\sum_{i=1}^{n} \alpha_i Var(X_i)$ for correlated random variables?

If there are no constraints on weights, there are no identical relationship among $Var(\sum_{i=1}^{n} \alpha_i X_i)$, $\sum_{i=1}^{n} \alpha_i^2 Var(X_i)$ and $\sum_{i=1}^{n} \alpha_i Var(X_i)$. For example, if all $\alpha_1, \cdots, \alpha_n$ are negative,

$$Var(\sum_{i=1}^{n} \alpha_i X_i) \ge 0 \ge \sum_{i=1}^{n} \alpha_i Var(X_i) .$$

And, various examples are as follows.

*Example 1* Suppose $X_1, \cdots, X_n$ are independent random variables,

If $\alpha_1 = \alpha_2 \cdots = \alpha_n = 1$, $Var(\sum_{i=1}^{n} \alpha_i X_i) = \sum_{i=1}^{n} \alpha_i Var(X_i)$.

If $\alpha_i > 1, i = 1, \cdots, n.$, $Var(\sum_{i=1}^{n} \alpha_i X_i) = \sum_{i=1}^{n} \alpha_i^2 Var(X_i) > \sum_{i=1}^{n} \alpha_i Var(X_i)$ .

If $0 \le \alpha_i \le 1, i = 1, \cdots, n.$, $Var(\sum_{i=1}^{n} \alpha_i X_i) = \sum_{i=1}^{n} \alpha_i^2 Var(X_i) < \sum_{i=1}^{n} \alpha_i Var(X_i)$ .

*Example 2.* Let $n = 2$ , $(X_1, X_2) \sim N(0,1;0,1;1)$ denote the 2-dimension Normal distribution.

If $\alpha_1 = \alpha_2 = 1$, $Var(\alpha_1 X_1 + \alpha_2 X_2) = 4 > \alpha_1 Var(X_1) + \alpha_2 Var(X_2) = 2$ .

If $\alpha_1=\alpha_2=\frac{1}{2}$, $Var(\alpha_1X_1+\alpha_2X_2)=1=\alpha_1Var(X_1)+\alpha_2Var(X_2)$.

If $\alpha_1=\frac{1}{2},\alpha_2=\frac{1}{3}$, $Var(\alpha_1X_1+\alpha_2X_2)=\frac{25}{36}<\alpha_1Var(X_1)+\alpha_2Var(X_2)=\frac{5}{6}$.

*Example 3.* Let $n=2$, $(X_1,X_2)\sim N(0,1;0,1;-1)$ denote the 2-dimension Normal distribution.

If $\alpha_1=\alpha_2=1$, $Var(\alpha_1X_1+\alpha_2X_2)=0<\alpha_1Var(X_1)+\alpha_2Var(X_2)=2$.

*Example 4.* . Let $n=2$, $\alpha_1=\alpha_2=\frac{1}{2}$, $(X_1,X_2)\sim N(0,1;0,1;1)$,

$$Var(\alpha_1X_1+\alpha_2X_2)=1>\sum_{i=1}^{2}\alpha_i^2Var(X_i)=\frac{1}{2}.$$

Let $n=2$, $\alpha_1=\alpha_2=\frac{1}{2}$, $(X_1,X_2)\sim N(0,1;0,1;-1)$,

$$Var(\alpha_1X_1+\alpha_2X_2)=0<\sum_{i=1}^{2}\alpha_i^2Var(X_i)=\frac{1}{2}.$$

From the above examples, we can conclude that, if all weights $\alpha_1,\cdots,\alpha_n$ are negative, or sum of all absolute values of weights $\alpha_1,\cdots,\alpha_n$ are larger than 1, the relationship among $Var(\sum_{i=1}^{n}\alpha_iX_i)$, $\sum_{i=1}^{n}\alpha_i^2Var(X_i)$ and $\sum_{i=1}^{n}\alpha_iVar(X_i)$ has no identical conclusion. Example 4 also illustrates that the lower bound of $Var(\sum_{i=1}^{n}\alpha_iX_i)\geq 0$ can be reached. Thus, we start our investigation from a simple case of non-negative weights $\alpha_1,\cdots,\alpha_n$ with sum of 1.

The rest of the paper is organized as follows. The variance inequalities are derived in Section 2. The Chebyshev's inequality of correlated random variables is obtained in Section 3. The law of large numbers of correlated random variables is obtained in Section 4. The conclusion is given in Section 5.

## 2. Variance inequality of correlated random variables

Consider the formula (3) , we first give the variance inequality with Cauchy-Schwarz inequality.

**Theorem 1.** Let $\alpha_1,\alpha_2,\cdots,\alpha_n$ be any real numbers with $0\leq\alpha_i\leq 1, i=1,\cdots,n$, and $\sum_{i=1}^{n}\alpha_i=1$.

Suppose $X_1,\cdots,X_n$ are any random variables with $Var(X_i)<+\infty$, $i=1,\cdots,n$, then

$$Var(\sum_{i=1}^{n}\alpha_i X_i) \le \sum_{i=1}^{n}\alpha_i Var(X_i). \qquad (5)$$

**Proof.** For any $X_i$ and $X_j$, we zoom the Cauchy-Schwarz inequality,

$$|Cov(X_i, X_j)| \le \sqrt{Var(X_i)}\sqrt{Var(X_j)} \le \frac{Var(X_i)+Var(X_j)}{2}.$$

Applying above formula to formula (3) and combine the same items, we have

$$\begin{aligned} Var(\sum_{i=1}^{n}\alpha_i X_i) &\le \sum_{i=1}^{n}\alpha_i^2 Var(X_i) + \sum_{1\le i<j\le n} 2\alpha_i\alpha_j \frac{Var(X_i)+Var(X_j)}{2} \\ &= \sum_{j=1}^{n}\alpha_j \sum_{i=1}^{n}\alpha_i Var(X_i) = (\sum_{i=1}^{n}\alpha_i)[\sum_{i=1}^{n}\alpha_i Var(X_i)] = \sum_{i=1}^{n}\alpha_i Var(X_i) \end{aligned} \qquad (6)$$

□

Note that using different Cauchy-Schwarz inequality form could obtain another upper bound of variance.

**Theorem 1'.** Let $\alpha_1, \alpha_2, \cdots, \alpha_n$ be any real numbers with $0 \le \alpha_i \le 1, i = 1, \cdots, n$, and $\sum_{i=1}^{n}\alpha_i = 1$.

Suppose $X_1, \cdots, X_n$ are any random variables with $Var(X_i) < +\infty$, $i = 1, \cdots, n$, then

$$Var(\sum_{i=1}^{n}\alpha_i X_i) \le (\sum_{i=1}^{n}\alpha_i^2)\sum_{i=1}^{n}Var(X_i). \qquad (5')$$

**Proof.** Denote $\sigma_i = \sqrt{Var(X_i)}$. According to Cauchy-Schwarz inequality

$$(\sum_{i=1}^{n}a_i b_i)^2 \le (\sum_{i=1}^{n}a_i^2)(\sum_{i=1}^{n}b_i^2),$$

$$\begin{aligned} Var(\sum_{i=1}^{n}\alpha_i X_i) &= \sum_{i=1}^{n}\sum_{j=1}^{n}\rho_{ij}\alpha_i\alpha_j\sqrt{Var(X_i)}\sqrt{Var(Y_i)} \le \sum_{i=1}^{n}\sum_{j=1}^{n}\alpha_i\alpha_j\sigma_i\sigma_j \\ &= (\sum_{i=1}^{n}\alpha_i\sigma_i)^2 \le (\sum_{i=1}^{n}\alpha_i^2)(\sum_{i=1}^{n}\sigma_i^2) = (\sum_{i=1}^{n}\alpha_i^2)(\sum_{i=1}^{n}Var(X_i)) \end{aligned}$$

□

To illustrate the difference between two upper bounds (5)(5'), we examine the values of $\sum_{i=1}^{n}\alpha_i Var(X_i)$ and $(\sum_{i=1}^{n}\alpha_i^2)\sum_{i=1}^{n}Var(X_i)$ in the range of $\alpha_i \in \{0.1, 0.2, \cdots, 1\}$, $\sum_{i=1}^{n}\alpha_i = 1$, $Var(X_i) \in \{0.1, 0.2, \cdots, 2\}$, $i = 1, \cdots, n$. The comparison results are list in Table 1.

As shown in **Table 1**, $\sum_{i=1}^{n}\alpha_i Var(X_i)$ is smaller than $(\sum_{i=1}^{n}\alpha_i^2)\sum_{i=1}^{n}Var(X_i)$ in most cases of simulations.

And, when $\alpha_1 = \alpha_2 = \cdots \alpha_n = \frac{1}{n}$, $\sum_{i=1}^{n} \alpha_i Var(X_i) = (\sum_{i=1}^{n} \alpha_i^2) \sum_{i=1}^{n} Var(X_i) = \frac{1}{n} \sum_{i=1}^{n} Var(X_i)$. For convenience, we discuss the upper bound (5) in the rest of the paper.

**Table 1**. Upper bound comparison simulation

| $n$ | Number of computation cases $\alpha_i \in \{0.1, 0.2, \cdots, 1\}$, $\sum_{i=1}^{n} \alpha_i = 1$ $Var(X_i) \in \{0.1, 0.2, \cdots, 2\}$, $i = 1, \cdots, n$ | Number of case $\sum_{i=1}^{n} \alpha_i Var(X_i) > (\sum_{i=1}^{n} \alpha_i^2) \sum_{i=1}^{n} Var(X_i)$ | Ratio (%) |
|---|---|---|---|
| 2 | 3600 | 520 | 14.44 |
| 3 | 288000 | 29137 | 10.11 |
| 4 | 13760000 | 799763 | 5.81 |

In fact, we can obtain the conclusion with matrix method. We give the proof in **Theorem 2.**

**Theorem 2.** Let $\alpha_1, \alpha_2, \cdots, \alpha_n$ be any real numbers, and $0 \le \alpha_i \le 1, i = 1, \cdots, n$, $\sum_{i=1}^{n} \alpha_i = 1$,

Suppose $X_1, \cdots, X_n$ are any random variables with $Var(X_i) < +\infty$, $i = 1, \cdots, n$, then

$$Var(\sum_{i=1}^{n} \alpha_i X_i) \le \sum_{i=1}^{n} \alpha_i Var(X_i). \tag{7}$$

**Proof.**

Denote $\sigma_i = \sqrt{Var(X_i)}$, $\rho_{ij} = \frac{Cov(X_i, X_j)}{\sqrt{Var(X_i)}\sqrt{Var(X_j)}}$, and $\rho_{ii} = 1$.

Then

$$Var(\sum_{i=1}^{n} \alpha_i X_i) - \sum_{i=1}^{n} \alpha_i Var(X_i) = \sum_{i=1}^{n} \sum_{i=1}^{n} Cov(\alpha_i X_i, \alpha_j X_j) - \sum_{i=1}^{n} \alpha_i Var(X_i)$$

$$= \sum_{i=1}^{n} \sum_{i=1}^{n} \rho_{ij} \alpha_i \alpha_j \sigma_i \sigma_j - \sum_{i=1}^{n} \alpha_i^2$$

$$= -(\sum_{i=1}^{n} \alpha_i \sigma_i^2 - \sum_{i=1}^{n} \sum_{i=1}^{n} \rho_{ij} \alpha_i \alpha_j \sigma_i \sigma_j)$$

We will prove that

$$\sum_{i=1}^{n} \alpha_i \sigma_i^2 - \sum_{i=1}^{n} \sum_{i=1}^{n} \rho_{ij} \alpha_i \alpha_j \sigma_i \sigma_j \ge 0. \tag{8}$$

Since $|\rho_{ij}| \le 1$, we have

$$\sum_{i=1}^{n}\alpha_i\sigma_i^2-\sum_{i=1}^{n}\sum_{i=1}^{n}\rho_{ij}\alpha_i\alpha_j\sigma_i\sigma_j \ge \sum_{i=1}^{n}\alpha_i\sigma_i^2-\sum_{i=1}^{n}\sum_{i=1}^{n}\alpha_i\alpha_j\sigma_i\sigma_j .$$

We will prove

$$\sum_{i=1}^{n}\alpha_i\sigma_i^2-\sum_{i=1}^{n}\sum_{i=1}^{n}\alpha_i\alpha_j\sigma_i\sigma_j \ge 0 \qquad (9).$$

In fact,

$$\sum_{i=1}^{n}\alpha_i\sigma_i^2-\sum_{i=1}^{n}\sum_{i=1}^{n}\alpha_i\alpha_j\sigma_i\sigma_j$$

$$=[\sigma_1,\cdots,\sigma_n]\begin{bmatrix}\alpha_1 & 0 & \cdots & 0\\ 0 & \alpha_2 & \cdots & 0\\ \cdots & \cdots & \cdots & \cdots\\ 0 & 0 & \cdots & \alpha_n\end{bmatrix}\begin{bmatrix}\sigma_1\\ \vdots\\ \sigma_n\end{bmatrix}-[\sigma_1,\cdots,\sigma_n]\begin{bmatrix}\alpha_1^2 & \alpha_1\alpha_2 & \cdots & \alpha_1\alpha_n\\ \alpha_1\alpha_2 & \alpha_2^2 & \cdots & \alpha_2\alpha_n\\ \cdots & \cdots & \cdots & \cdots\\ \alpha_1\alpha_n & \alpha_2\alpha_n & \cdots & \alpha_n^2\end{bmatrix}\begin{bmatrix}\sigma_1\\ \vdots\\ \sigma_n\end{bmatrix}$$

$$=[\sigma_1,\cdots,\sigma_n]\left(\begin{bmatrix}\alpha_1 & 0 & \cdots & 0\\ 0 & \alpha_2 & \cdots & 0\\ \cdots & \cdots & \cdots & \cdots\\ 0 & 0 & \cdots & \alpha_n\end{bmatrix}-\begin{bmatrix}\alpha_1^2 & \alpha_1\alpha_2 & \cdots & \alpha_1\alpha_n\\ \alpha_1\alpha_2 & \alpha_2^2 & \cdots & \alpha_2\alpha_n\\ \cdots & \cdots & \cdots & \cdots\\ \alpha_1\alpha_n & \alpha_2\alpha_n & \cdots & \alpha_n^2\end{bmatrix}\right)\begin{bmatrix}\sigma_1\\ \vdots\\ \sigma_n\end{bmatrix}$$

$$=[\sigma_1,\cdots,\sigma_n]\begin{bmatrix}\alpha_1-\alpha_1^2 & -\alpha_1\alpha_2 & \cdots & -\alpha_1\alpha_n\\ -\alpha_1\alpha_2 & \alpha_2-\alpha_2^2 & \cdots & -\alpha_2\alpha_n\\ \cdots & \cdots & \cdots & \cdots\\ -\alpha_1\alpha_n & -\alpha_2\alpha_n & \cdots & \alpha_n-\alpha_n^2\end{bmatrix}\begin{bmatrix}\sigma_1\\ \vdots\\ \sigma_n\end{bmatrix}$$

This is a quadratic form. We only prove that Matrix

$$A=\begin{bmatrix}\alpha_1-\alpha_1^2 & -\alpha_1\alpha_2 & \cdots & -\alpha_1\alpha_n\\ -\alpha_1\alpha_2 & \alpha_2-\alpha_2^2 & \cdots & -\alpha_2\alpha_n\\ \cdots & \cdots & \cdots & \cdots\\ -\alpha_1\alpha_n & -\alpha_2\alpha_n & \cdots & \alpha_n-\alpha_n^2\end{bmatrix}$$

is a positive semidefinite matrix. Consider the $k$ -th ($1\le k\le n$) principal minors.

$$A\begin{pmatrix}i_1 & i_2 & \cdots & i_k\\ i_1 & i_2 & \cdots & i_k\end{pmatrix}=\begin{vmatrix}\alpha_{i_1}-\alpha_{i_1}^2 & -\alpha_{i_1}\alpha_{i_2} & \cdots & -\alpha_{i_1}\alpha_{i_k}\\ -\alpha_{i_1}\alpha_{i_2} & \alpha_{i_2}-\alpha_{i_2}^2 & \cdots & -\alpha_{i_2}\alpha_{i_k}\\ \cdots & \cdots & \cdots & \cdots\\ -\alpha_{i_1}\alpha_k & -\alpha_{i_2}\alpha_{i_k} & \cdots & \alpha_{i_k}-\alpha_{i_k}^2\end{vmatrix}=\prod_{m=1}^{k}\alpha_{i_m}\begin{vmatrix}1-\alpha_{i_1} & -\alpha_{i_1} & \cdots & -\alpha_{i_1}\\ -\alpha_{i_2} & 1-\alpha_{i_2} & \cdots & -\alpha_{i_2}\\ \cdots & \cdots & \cdots & \cdots\\ -\alpha_{i_k} & -\alpha_{i_k} & \cdots & 1-\alpha_{i_k}\end{vmatrix}$$

$$=(\prod_{m=1}^{k}\alpha_{i_m})\begin{vmatrix}1-\sum_{m=1}^{k}\alpha_{i_m} & 1-\sum_{m=1}^{k}\alpha_{i_m} & \cdots & 1-\sum_{m=1}^{k}\alpha_{i_m}\\ -\alpha_{i_2} & 1-\alpha_{i_2} & \cdots & -\alpha_{i_2}\\ \cdots & \cdots & \cdots & \cdots\\ -\alpha_{i_k} & -\alpha_{i_k} & \cdots & 1-\alpha_{i_k}\end{vmatrix}=(\prod_{m=1}^{k}\alpha_{i_m})(1-\sum_{m=1}^{k}\alpha_{i_m})\begin{vmatrix}1 & 1 & \cdots & 1\\ -\alpha_{i_2} & 1-\alpha_{i_2} & \cdots & -\alpha_{i_2}\\ \cdots & \cdots & \cdots & \cdots\\ -\alpha_{i_k} & -\alpha_{i_k} & \cdots & 1-\alpha_{i_k}\end{vmatrix}$$

$$=(\prod_{m=1}^{k}\alpha_{i_m})(1-\sum_{m=1}^{k}\alpha_{i_m})\begin{vmatrix}1 & 1 & \cdots & 1\\ 0 & 1 & \cdots & 0\\ \cdots & \cdots & \cdots & \cdots\\ 0 & 0 & \cdots & 1\end{vmatrix}=(\prod_{m=1}^{k}\alpha_{i_m})(1-\sum_{m=1}^{k}\alpha_{i_m})\geq 0 \tag{10}$$

Then, $A$ is positive semidefinite. Thus, formula (8) holds, and formula (7) holds. This ends of proof . □

If $0\leq\alpha_i\leq 1, i=1,\cdots,n, 0<\sum_{i=1}^{n}\alpha_i<1,$ we have the following theorem.

**Theorem 3.** Let $\alpha_1,\alpha_2,\cdots,\alpha_n$ be real numbers, $0\leq\alpha_i\leq 1, i=1,\cdots,n, 0\leq\sum_{i=1}^{n}\alpha_i<1$ . Suppose $X_1,\cdots,X_n$ are any random variables with $Var(X_i)<+\infty$, $i=1,\cdots,n$ then

$$Var(\sum_{i=1}^{n}\alpha_i X_i)\leq\sum_{i=1}^{n}\alpha_i Var(X_i). \tag{11}$$

**Proof.** Examine the formula (10) in the proof of Theorem 1, the condition $0\leq\alpha_i\leq 1$ , $0<\sum_{i=1}^{n}\alpha_i<1$, can guarantee that conclusion of formula (8) is still true. Hence, Theorem 3 holds.

□

While $\alpha_i\geq 0, i=1,\cdots,n, \sum_{i=1}^{n}\alpha_i>0$, since the weights can be normalized to sum 1, the upper bound of $Var(\sum_{i=1}^{n}\alpha_i X_i)$ is described as follows.

**Theorem 4** Let $\alpha_1,\alpha_2,\cdots,\alpha_n$ be real numbers, $\alpha_i\geq 0, i=1,\cdots,n$ . Suppose $X_1,\cdots,X_n$ are any random variables with $Var(X_i)<+\infty$, $i=1,\cdots,n$, then

$$Var(\sum_{i=1}^{n}\alpha_i X_i) \le (\sum_{i=1}^{n}\alpha_i)\sum_{i=1}^{n}\alpha_i Var(X_i) \quad (12)$$

**Proof.** If $\alpha_i = 0, i = 1,\cdots,n$, formula (12) holds.

If $\sum_{i=1}^{n}\alpha_i > 0$, according to Theorem 2, we obtain

$$Var(\sum_{i=1}^{n}\alpha_i X_i) = Var[(\sum_{i=1}^{n}\alpha_i)\sum_{i=1}^{n}\frac{\alpha_i}{\sum_{i=1}^{n}\alpha_i}X_i] = (\sum_{i=1}^{n}\alpha_i)^2 Var[\sum_{i=1}^{n}\frac{\alpha_i}{\sum_{i=1}^{n}\alpha_i}X_i]$$
$$\le (\sum_{i=1}^{n}\alpha_i)^2\sum_{i=1}^{n}\frac{\alpha_i}{\sum_{i=1}^{n}\alpha_i}Var(X_i) \le (\sum_{i=1}^{n}\alpha_i)\sum_{i=1}^{n}\alpha_i Var(X_i)$$

□

In fact, the conclusion of Theorem 4 can also be obtained by checking the formula (6) in Theorem 1. Furthermore, since $Var(\pm X) = Var(X)$ and,

$$\sum_{i=1}^{n}\alpha_i X_i = \sum_{i=1}^{n}|\alpha_i|[\mathrm{sgn}(\alpha_i)X_i], \quad (13)$$

where $\mathrm{sgn}(\bullet)$ is the sign function. We extend the weights to general case without any limitation, and have the following conclusion.

**Theorem 5** Let $\alpha_1,\cdots,\alpha_n$ be any $n$ real numbers. Suppose $X_1,\cdots,X_n$ are any random variables with $Var(X_i) < +\infty$, $i = 1,\cdots,n$, then

$$Var(\sum_{i=1}^{n}\alpha_i X_i) \le (\sum_{i=1}^{n}|\alpha_i|)\sum_{i=1}^{n}|\alpha_i| Var(X_i) \quad (14)$$

**Proof.** If $\alpha_i = 0, i = 1,\cdots,n$, formula (14) holds.

If $\sum_{i=1}^{n}|\alpha_i| > 0$, according to Theorem 4, we obtain

$$Var(\sum_{i=1}^{n}\alpha_i X_i) = \sum_{i=1}^{n}Var(|\alpha_i|[\mathrm{sgn}(\alpha_i)X_i]) \le (\sum_{i=1}^{n}|\alpha_i|)[\sum_{i=1}^{n}|\alpha_i| Var(\mathrm{sgn}(\alpha_i)X_i)]$$
$$= (\sum_{i=1}^{n}|\alpha_i|)[\sum_{i=1}^{n}|\alpha_i| Var(X_i)]$$

□

Obviously, we can obtain the following corollaries.

**Corollary 1.** Let $X_1,\cdots,X_n$ be any $n$ random variables. Suppose $X_1,\cdots,X_n$ are any random variables with $Var(X_i)<+\infty$, $i=1,\cdots,n$, then

$$Var(\frac{1}{n}\sum_{i=1}^{n}X_i)\le\frac{1}{n}\sum_{i=1}^{n}Var(X_i) \tag{15}$$

**Proof** Let $\alpha_i=\frac{1}{n}, i=1,\cdots,n$, according to Theorem 2, the conclusion holds. □

Corollary 1 shows that variance of mean could not exceed mean of variances for correlated random variables, it is also true for independent random variables through checking formula (4) and Example 1.

**Corollary 2 If** $0\le|\alpha_i|\le 1, i=1,\cdots,n, \sum_{i=1}^{n}|\alpha_i|\le 1,$ for any random variables $X_1,\cdots,X_n$ with $Var(X_i)<+\infty$, $i=1,\cdots,n$, we have

$$Var(\sum_{i=1}^{n}\alpha_i X_i)\le(\sum_{i=1}^{n}|\alpha_i|)\sum_{i=1}^{n}|\alpha_i|Var(X_i)\le\sum_{i=1}^{n}|\alpha_i|Var(X_i)\le\sum_{i=1}^{n}Var(X_i) \tag{16}$$

**Proof.** As $0\le|\alpha_i|\le 1, i=1,\cdots,n, \sum_{i=1}^{n}|\alpha_i|\le 1,$ according to **Theorem 5,** the inequalities are true from left to right. □

If the variables are uncorrelated (or independent), the relationship among $Var(\sum_{i=1}^{n}\alpha_i X_i)$, $\sum_{i=1}^{n}\alpha_i^2 Var(X_i)$, and $\sum_{i=1}^{n}\alpha_i Var(X_i)$ are as follows.

**Corollary 3 If** $0\le|\alpha_i|\le 1, i=1,\cdots,n, \sum_{i=1}^{n}|\alpha_i|\le 1,$ and $X_1,\cdots,X_n$ be any $n$ uncorrelated ( or independent) random variables,

$$Var(\sum_{i=1}^{n}\alpha_i X_i)=\sum_{i=1}^{n}\alpha_i^2 Var(X_i)\le\sum_{i=1}^{n}|\alpha_i|Var(X_i)\le\sum_{i=1}^{n}Var(X_i) \tag{17}$$

**Proof.** According to formula (3) and $\alpha_i^2\le|\alpha_i|$, while $0\le|\alpha_i|\le 1$, $i=1,\cdots,n$. The conclusion is obviously true. □

In addition, an inequality for the covariance and variance among correlated random variables is achieved as follows.

**Corollary 4.** Let $X_1, \cdots, X_n$ be any $n$ correlated random variables, with $Var(X_i) < +\infty$, $i = 1, \cdots, n$. Then

$$-\frac{1}{n}\sum_{i=1}^{n} Var(X_i) \le \frac{2}{n}\sum_{1\le i<j\le n} Cov(X_i, X_j) \le (1-\frac{1}{n})\sum_{i=1}^{n} Var(X_i) \tag{18}$$

**Proof.**

$$Var(\frac{1}{n}\sum_{i=1}^{n} X_i) = \frac{1}{n^2} Var(\sum_{i=1}^{n} X_i) = \frac{1}{n^2}[\sum_{i=1}^{n} Var(X_i) + 2\sum_{1\le i<j\le n} Cov(X_i, X_j)]$$

$$= \frac{1}{n}[\frac{1}{n}\sum_{i=1}^{n} Var(X_i) + \frac{2}{n}\sum_{1\le i<j\le n} Cov(X_i, X_j)] \tag{19}$$

As $Var(\frac{1}{n}\sum_{i=1}^{n} X_i) \ge 0$, the left inequality is proved.

Using formula (15) $Var(\frac{1}{n}\sum_{i=1}^{n} X_i) \le \frac{1}{n} Var(\sum_{i=1}^{n} X_i)$ and formula (19), the right inequality is obtained.

□

# 3. Chebyshev inequality of correlated random variables

Chebyshev's inequality is an important probability formula, when sum of correlated random variables are considered, we have the following inequality.

**Theorem 6** (Chebyshev inequality) Let $\alpha_1, \cdots, \alpha_n$ be any $n$ real numbers, and $X_1, \cdots, X_n$ be any $n$ correlated random variables with $Var(X_i) < +\infty$, $i = 1, \cdots, n$. For any $\delta > 0$, we have

$$P(|\sum_{i=1}^{n} \alpha_i X_i - \sum_{i=1}^{n} \alpha_i EX_i| > \delta) \le \frac{(\sum_{i=1}^{n} |\alpha_i|)}{\delta^2}\sum_{i=1}^{n} |\alpha_i| Var(X_i) \tag{20}$$

**Proof.**

According to Chebyshev inequality and Theorem 5, we obtain

$$P(|\sum_{i=1}^{n} \alpha_i X_i - \sum_{i=1}^{n} \alpha_i EX_i| > \delta) \le \frac{Var(\sum_{i=1}^{n} \alpha_i X_i)}{\delta^2} \le \frac{(\sum_{i=1}^{n} |\alpha_i|)}{\delta^2}\sum_{i=1}^{n} |\alpha_i| Var(X_i)$$

□

Let $\alpha_i = \frac{1}{n}, i = 1, \cdots, n$, according to formula (20), we have the following conclusion.

**Corollary 5.** Let $X_1, \cdots, X_n$ be any $n$ correlated random variables with $Var(X_i) < +\infty$, $i = 1, \cdots, n$. For any $\delta > 0$,

$$P(|\frac{1}{n}\sum_{i=1}^{n} X_i - \frac{1}{n}\sum_{i=1}^{n} EX_i| > \delta) \le \frac{1}{\delta^2}[\frac{1}{n}\sum_{i=1}^{n} Var(X_i)] \quad (21)$$

When $X_1, \cdots, X_n$ are uncorrelated (or independent) random variables, the Chebyshev's inequality is

$$P(|\frac{1}{n}\sum_{i=1}^{n} X_i - \frac{1}{n}\sum_{i=1}^{n} EX_i| > \delta) \le \frac{1}{\delta^2}[\frac{1}{n^2}\sum_{i=1}^{n} Var(X_i)] \quad (22)$$

Formula (21) illustrates a relationship among $\frac{1}{n}\sum_{i=1}^{n} X_i$, $\frac{1}{n}\sum_{i=1}^{n} EX_i$ and $\frac{1}{n}\sum_{i=1}^{n} Var(X_i)$, and it could be applied to probability estimation of $\frac{1}{n}\sum_{i=1}^{n} X_i$, regardless of the correlation among $X_1, \cdots, X_n$. When $X_1, \cdots, X_n$ are uncorrelated or independent, the formula (21) sharpens to formula (22). This conclusion also demonstrates that correlation among variables enlarges the departure of $\frac{1}{n}\sum_{i=1}^{n} X_i$ from its expectation.

**Corollary 6.** Let $X_1, \cdots, X_n$ be any $n$ correlated random variables, with $0 < Var(X_i) < +\infty$, $i = 1, \cdots, n$. For any $\delta > 0$,

$$P(|\frac{1}{n}\sum_{i=1}^{n} \frac{X_i - EX_i}{\sqrt{VarX_i}}| > \delta) \le \frac{1}{\delta^2} \quad (23)$$

$$P(|\sum_{i=1}^{n} \frac{X_i - EX_i}{\sqrt{VarX_i}}| > \delta) \le \frac{n^2}{\delta^2} \quad (24)$$

**Proof.** According to formula (21)

$$P(|\frac{1}{n}\sum_{i=1}^{n} \frac{X_i - EX_i}{\sqrt{VarX_i}}| > \delta) \le \frac{1}{\delta^2}[\frac{1}{n}\sum_{i=1}^{n} Var(\frac{X_i}{\sqrt{VarX_i}})] = \frac{1}{\delta^2}.$$

Thus, formula (23) holds.

According to formula (20), for variables $\frac{X_1}{\sqrt{VarX_1}},\cdots,\frac{X_n}{\sqrt{VarX_n}}$ with weights $\alpha_1=\alpha_2=\cdots=\alpha_n=1$, we have

$$P(|\sum_{i=1}^{n}\frac{X_i-EX_i}{\sqrt{VarX_i}}|>\delta)\le\frac{1}{\delta^2}[n\sum_{i=1}^{n}Var(\frac{X_i}{\sqrt{VarX_i}})]=\frac{n^2}{\delta^2}.$$

Then, formula (24) holds.

This ends of the proof. □

Obviously, an alternative proof of formula (24) is easily derived from formula (23) as follows.

$$P(|\sum_{i=1}^{n}\frac{X_i-EX_i}{\sqrt{VarX_i}}|>\delta)=P(|\frac{1}{n}\sum_{i=1}^{n}\frac{X_i-EX_i}{\sqrt{VarX_i}}|>\frac{\delta}{n})\le\frac{1}{(\frac{\delta}{n})^2}\le\frac{n^2}{\delta^2}.$$

The Chebyshev inequalities in **Theorem 6, Corollary** 5 and **Corollary** 6 are different from the Chebyshev inequality in Theorem 1.1 in [7]. In [7], the Chebyshev inequality models the probability of "$(X-EX)^T\Sigma^{-1}(X-EX)>\varepsilon$", where $X=(X_1,\cdots,X_n)^T$, $\Sigma$ is a positive matrix, its inverse event "$(X-EX)^T\Sigma^{-1}(X-EX)\le\varepsilon$" can be treated as the random event of the vector $X$ from expectation $EX$ with a Mahalanobis distance in high dimension random vector space [8]. And, "$(X-EX)^T\Sigma^{-1}(X-EX)\le\varepsilon$" means an $n$-dimension ellipsoid neighborhood, while the event "$|\sum_{i=1}^{n}\frac{X_i-EX_i}{\sqrt{VarX_i}}|\le\delta$"represents a more complex neighborhood. Since, the event "$\sum_{i=1}^{n}|\frac{X_i-EX_i}{\sqrt{VarX_i}}|\le\delta$" represents an $n$-dimension cube neighborhood of scaled random vector $(\frac{X_1}{\sqrt{VarX_1}},\frac{X_2}{\sqrt{VarX_2}},\cdots,\frac{X_n}{\sqrt{VarX_n}})$ from expectation $(\frac{EX_1}{\sqrt{VarX_1}},\frac{EX_2}{\sqrt{VarX_2}},\cdots,\frac{EX_n}{\sqrt{VarX_n}})$, and "$\sum_{i=1}^{n}|\frac{X_i-EX_i}{\sqrt{VarX_i}}|\le\delta$" $\subset$ "$|\sum_{i=1}^{n}\frac{X_i-EX_i}{\sqrt{VarX_i}}|\le\delta$", the event "$|\sum_{i=1}^{n}\frac{X_i-EX_i}{\sqrt{VarX_i}}|>\delta$" is included in the domain of event "$\sum_{i=1}^{n}|\frac{X_i-EX_i}{\sqrt{VarX_i}}|>\delta$". The intuitionistic domain of "$|\sum_{i=1}^{n}\frac{X_i-EX_i}{\sqrt{VarX_i}}|>\delta$" in this paper is slight complex. Hence, **Theorem 6, Corollary** 5 and **Corollary** 6 address another Chebyshev inequality of high dimensional random variables.

In other words, for the sequence of random variables $X_1,\cdots,X_n,\cdots$, let $X=(X_1,\cdots,X_n)^T$, $S_n=\frac{1}{n}\sum_{i=1}^{n}X_i$ . $(X-EX)^T\Sigma^{-1}(X-EX)$ is the quadratic form of $X_1,\cdots,X_n$, it's a $n$-dimension distance. However, $|S_n-ES_n|$ is the one-dimension distance of sum $S_n$ from its expectation $ES_n$. The representation difference of mathematical formulae also shows the different between Theorem 6, Corollary 5 and Corollary 6 and Theorem 1.1 in [7].

## 4. Law of large numbers of correlated random variables

Suppose $X_1,\cdots,X_n,\cdots$ are any random variable sequences, let $S_n=\frac{1}{n}\sum_{i=1}^{n}X_i$, for any $\varepsilon>0$, if

$$\frac{1}{n^2}Var(\sum_{i=1}^{n}X_i)\to 0 \qquad (n\to\infty), \tag{25}$$

Then,

$$\lim_{n\to\infty}P(|\frac{1}{n}\sum_{i=1}^{n}X_i-\frac{1}{n}\sum_{i=1}^{n}EX_i|<\varepsilon)=1 \tag{26}$$

(25) is the Markov sufficient condition [6,9] of (26). (26) is denoted as $\frac{1}{n}\sum_{i=1}^{n}X_i\xrightarrow{P}\frac{1}{n}\sum_{i=1}^{n}EX_i$, and called weak law of large numbers (WLLN).

If any of two variables of $X_1,\cdots,X_n,\cdots$ are uncorrelated, and there exists a constant $C>0$,

$$DX_i\le C, \tag{27}$$

(26) holds. It is called Chebyshev Law of large numbers [6,9].

When $X_1,\cdots,X_n,\cdots$ are correlated, the bound condition (27) could not ensure that (25) is true, especially when the correlations are unknown. According to (15), we have a Markov Law of large numbers.

**Theorem 7** (Markov Law of large numbers) Let $X_1,\cdots,X_n,\cdots$ be any correlated random variables sequences with $Var(X_i)<+\infty$, $i=1,2,\cdots$. For any $\varepsilon>0$, if

$$\frac{1}{n}\sum_{i=1}^{n} Var(X_i) \to 0 \qquad (n \to +\infty) \tag{28}$$

Then,

$$\lim_{n\to\infty} P(|\frac{1}{n}\sum_{i=1}^{n} X_i - \frac{1}{n}\sum_{i=1}^{n} EX_i| < \varepsilon) = 1$$

**Proof.** According to (15)

$$\frac{1}{n^2} Var(\sum_{i=1}^{n} X_i) = Var(\frac{1}{n}\sum_{i=1}^{n} X_i) \le \frac{1}{n}\sum_{i=1}^{n} Var(X_i)$$

Using (28) and (21), we have $\lim_{n\to\infty} P(|\frac{1}{n}\sum_{i=1}^{n} X_i - \frac{1}{n}\sum_{i=1}^{n} EX_i| < \varepsilon) = 1$.

□

(28) is a sufficient condition for law of large numbers without known the correlation between variables. Theorem 7 provides a convenient condition to verify the law of large numbers or probability convergence of correlated random variable sequence.

Furthermore, we can provide another sufficient condition with the following Lemma.

**Lemma 1.** Let $X_1, \cdots, X_n, \cdots$ be any correlated random variables sequence with $Var(X_i) < +\infty$, $i = 1, 2, \cdots$. For any $0 < r \le s$, we have,

$$[\frac{1}{n}\sum_{i=1}^{n} (Var(X_i))^r]^{\frac{1}{r}} \le [\frac{1}{n}\sum_{i=1}^{n} (Var(X_i))^s]^{\frac{1}{s}} \tag{29}$$

**Proof.** For nonnegative values $Var(X_1), \cdots, Var(X_n)$, we define a random variable $\xi$ with probability distribution

$$P(\xi = Var(X_i)) = \frac{1}{n}, \qquad i = 1, \cdots, n.$$

For any $0 < r \le s$, using the Lyapunov's inequality [3], we have

$$(E|\xi|^r)^{\frac{1}{r}} \le (E|\xi|^s)^{\frac{1}{s}},$$

which is

$$[\sum_{i=1}^{n} \frac{1}{n} (Var(X_i))^r]^{\frac{1}{r}} \le [\sum_{i=1}^{n} \frac{1}{n} (Var(X_i))^s]^{\frac{1}{s}}.$$

Hence, (29) holds. □

Utilizing (29), we can have the following theorem.

**Theorem 8** (Markov Law of large numbers) Let $X_1,\cdots,X_n,\cdots$ be any correlated random variables sequences with $Var(X_i)<+\infty$, $i=1,2,\cdots$. For any $\varepsilon>0$, $s\geq 1$, if

$$[\frac{1}{n}\sum_{i=1}^{n}(Var(X_i))^s]^{\frac{1}{s}}\to 0 \qquad (n\to+\infty) \tag{30}$$

Then,

$$\lim_{n\to\infty}P(|\frac{1}{n}\sum_{i=1}^{n}X_i-\frac{1}{n}\sum_{i=1}^{n}EX_i|<\varepsilon)=1$$

**Proof.** Let $r=1$, $s\geq 1$ in (27), we have

$$\frac{1}{n}\sum_{i=1}^{n}Var(X_i)\leq[\frac{1}{n}\sum_{i=1}^{n}(Var(X_i))^s]^{\frac{1}{s}} \tag{31}$$

According to (30)(31), the condition (25) holds. Applying Theorem 7, we have

$$\lim_{n\to\infty}P(|\frac{1}{n}\sum_{i=1}^{n}X_i-\frac{1}{n}\sum_{i=1}^{n}EX_i|<\varepsilon)=1.$$

□

Now, we give an example of correlated random variables sequence.

***Example 5.*** Let $X_1,\cdots,X_n,\cdots$ be independent and identical distribution $N(\mu,\sigma^2)$ ($\sigma>0$).

Let $S_n=\frac{1}{n}\sum_{i=1}^{n}X_i$, then $E(S_n)=\mu, Var(S_n)=\frac{\sigma^2}{n}$, $n=1,2,\cdots$.

$$\begin{aligned}
Cov(S_i,S_{i+k}) &= E(S_i-ES_i)(S_{i+k}-ES_{i+k})=E(S_iS_{i+k})-ES_iES_{i+k}\\
&=E[\frac{1}{i}\sum_{j=1}^{i}X_j(\frac{1}{i+k}\sum_{j=1}^{i}X_j+\frac{1}{i+k}\sum_{j=i+1}^{i+k}X_j)]-ES_iES_{i+k}\\
&=\frac{1}{i}\frac{1}{i+k}E(\sum_{j=1}^{i}X_j)^2+\frac{1}{i}\frac{1}{i+k}E(\sum_{j=1}^{i}X_j\sum_{j=i+1}^{i+k}X_j)-ES_iES_{i+k}\\
&=\frac{\sigma^2+i\mu^2}{(i+k)}+\frac{1}{i}\frac{1}{i+k}ik\mu^2-\mu^2=\frac{\sigma^2}{(i+k)}
\end{aligned}$$

For the correlated random variable sequence $S_1,S_2,\cdots,S_n,\cdots$, let $Y_n=\frac{1}{n}\sum_{i=1}^{n}S_i$. We have,

$$EY_n=E(\frac{1}{n}\sum_{i=1}^{n}S_i)=\mu,$$

The Markov sufficient condition (25),

$$Var(Y_n)=Var(\frac{1}{n}\sum_{i=1}^{n}S_i)=\frac{1}{n^2}Var(\sum_{i=1}^{n}S_i)=\frac{1}{n^2}\sum_{i=1}^{n}\sum_{j=1}^{n}Cov(S_i,S_j)=\frac{1}{n^2}[\sum_{i=1}^{n}\frac{\sigma^2}{i}+2\sum_{1\le i<j\le n}\frac{\sigma^2}{j}]$$
$$=\frac{1}{n^2}[\sum_{i=1}^{n}\frac{\sigma^2}{i}+2\sum_{i=1}^{n-1}\sum_{j=i+1}^{n}\frac{\sigma^2}{j}]=\frac{1}{n^2}[\sum_{i=1}^{n}\frac{\sigma^2}{i}+2\sum_{j=2}^{n}(j-1)\frac{\sigma^2}{j}]$$
$$=\frac{1}{n^2}[\sum_{i=1}^{n}\frac{\sigma^2}{i}+2(n-1)\sigma^2-2\sum_{j=2}^{n}\frac{\sigma^2}{j}]=\frac{1}{n^2}[2n\sigma^2-\sum_{i=1}^{n}\frac{\sigma^2}{i}]$$
$$=\frac{\sigma^2}{n}[2-\frac{1}{n}(C+\ln n+\varepsilon_n)]\to 0,\quad n\to+\infty$$

The Markov sufficient condition (28) in Theorem 7,

$$Var(Y_n)\le\frac{1}{n}\sum_{i=1}^{n}Var(S_i)=\frac{1}{n}\sum_{i=1}^{n}\frac{\sigma^2}{i}=\frac{\sigma^2}{n}(C+\ln n+\varepsilon_n)\to 0,\quad n\to+\infty$$

where $C=0.577216\cdots$ is Euler constant, $\varepsilon_n\to 0,\quad n\to+\infty$.

Then, for any $\varepsilon>0$, $\lim_{n\to\infty}P(|\frac{1}{n}\sum_{i=1}^{n}S_i-\frac{1}{n}\sum_{i=1}^{n}ES_i|<\varepsilon)=1$.

Example 5 shows that the law of large numbers of the correlated random variable sequence $S_1,S_2,\cdots,S_n,\cdots$ holds, and the sufficient condition (28) is simpler to be calculated than the sufficient condition (25). Note that Theorem 7 works in the case of non-stationary time series [11,12]. For the stationary time series, Theorem 7 fails.

As for the correlated random variables satisfying (27), we have the following theorem.

**Theorem 9** (Chebyshev Law of large numbers) Let $X_1,\cdots,X_n,\cdots$ be any correlated random variables sequence. If there exists a constant $C>0$, $Var(X_i)<C$, $i=1,2,\cdots$, for any $\varepsilon>0$, $s>0$, then,

$$\lim_{n\to\infty}P(|\frac{1}{n^{1+s}}\sum_{i=1}^{n}X_i-\frac{1}{n^{1+s}}\sum_{i=1}^{n}EX_i|<\varepsilon)=1$$

**Proof:**

$$Var(\frac{1}{n^{1+s}}\sum_{i=1}^{n}X_i)=\frac{1}{n^{2+2s}}Var(\sum_{i=1}^{n}X_i)\le\frac{1}{n^{2+2s}}n\sum_{i=1}^{n}Var(X_i)\le\frac{1}{n^{1+2s}}\sum_{i=1}^{n}Var(X_i)\quad (32)$$
$$\le\frac{1}{n^{1+2s}}nC=\frac{1}{n^{2s}}C\to 0,\quad n\to\infty$$

According to (21)(32), we have

$$\lim_{n\to\infty} P\left(\left|\frac{1}{n^{1+s}}\sum_{i=1}^{n} X_i - \frac{1}{n^{1+s}}\sum_{i=1}^{n} EX_i\right| < \varepsilon\right) = 1.$$

□

As a special case of Theorem 9, we can obtain the theorem for the random variables with same expectation and variance.

**Theorem 10** Let $X_1, \cdots, X_n, \cdots$ be any correlated random variable sequence with same expectation $E(X_i) = \mu$ and variance $Var(X_i) = \sigma^2$, $i = 1, 2, \cdots$, for any $\varepsilon > 0$, $s > 0$, then

$$\lim_{n\to\infty} P\left(\left|\frac{1}{n^{1+s}}\sum_{i=1}^{n} X_i - \frac{1}{n^{1+s}}\sum_{i=1}^{n} EX_i\right| < \varepsilon\right) = 1.$$

**Proof.** The proof is similar to Theorem 9. We omitted it. □

Corresponding to the famous Bernoulli Law of large numbers of independent and identical distribution (iid) random varables , we have

**Theorem 11**(Bernoulli Law of large numbers) Let $X_1, \cdots, X_n, \cdots$ be any correlated random variables sequence with identically distribution:

$P(X_i = 1) = p$, $P(X_i = 0) = 1 - p$, $i = 1, 2, \cdots$, $(0 < p < 1)$

For any $\varepsilon > 0$, $s > 0$, then

$$\lim_{n\to\infty} P\left(\left|\frac{1}{n^{1+s}}\sum_{i=1}^{n} X_i - \frac{p}{n^s}\right| < \varepsilon\right) = 1.$$

**Theorem** 9,10,11 discuss the stable of $\{\frac{1}{n^{1+s}}\sum_{i=1}^{n} X_i\}$, motivated by Bernoulli Law of large numbers in 1713, we are still interested in the stable of $\{\frac{1}{n}\sum_{i=1}^{n} X_i\}$, as its significance lies in not only the extension of Bernoulli WLLN, but also the moment method of estimation of point estimation [13], we continue to discuss the WLLN of correlated random variables.

To illustrate our motivation, we first give an example of correlated random variables of 0-1 distribution $b(1, p)$.

***Example 6.***（Random telegraph signals ([14],p24)） Let $\{X(t), t\in(-\infty,+\infty)\}$ be a stochastic process with state space $S=\{0,1\}$， which satisfies

(1) For any $t\in(-\infty,+\infty)$, $P(X(t)=0)=P(X(t)=1)=\frac{1}{2}$.

(2) For any $t$, the number of zero crossings in the interval $(\mathrm{t},\ \mathrm{t}+\Delta\mathrm{t})$ is described by a Poisson process $N(t+\Delta t)-N(t)\sim\Pi(\lambda\Delta t)$, $\lambda>0$.

(3) $X(t)$ is independent to $N(t+\Delta t)-N(t)$.

$\{X(t), t\in(-\infty,+\infty)\}$ are called random telegraph signals. The special case $\{X(n), n=0,1,2,\cdots\cdots\}$ can be treated as the correlated case of Bernoulli sequences. We discuss the WLLN of $\{X(n), n=0,1,2,\cdots\cdots\}$.

The expected value is

$$E[X(t)]=0\times P\{X(t)=0\}+1\times P\{X(t)=1\}=\frac{1}{2}$$

For any $\tau\ge 0$,

$$\{X(t)X(t+\tau)=1\}=\{X(t)=1, X(t+\tau)=1\}=\{X(t)=1, N(t+\tau)-N(t)=\text{even number}\}$$

$$=\{X(t)=1, \bigcup_{k=0}^{\infty}\{N(\tau)=2k\}\}$$

$$P\{X(t)X(t+\tau)=1\}=P\{X(t)=1\}P\left\{\bigcup_{k=0}^{\infty}N(\tau)=2k\right\}=\frac{1}{2}\sum_{k=0}^{\infty}P\{N(\tau)=2k\}=\frac{1}{4}[1+e^{-2\lambda\tau}]$$

The autocorrelation function is

$$R_X(t,t+\tau)=E[X(t)X(t+\tau)]=1\times P\{X(t)X(t+\tau)=1\}+0\times P\{X(t)X(t+\tau)=0\}=\frac{1}{4}[1+e^{-2\lambda\tau}]$$

$$Cov(X(t), X(t+\tau))=\frac{1}{4}e^{-2\lambda\tau}$$

Therefore, {X(n)} is correlated sequence.

$$D(\frac{1}{n}\sum_{k=1}^{n}X_k)=\frac{1}{n^2}\sum_{i=1}^{n}\sum_{j=1}^{n}Cov(X_i,X_j)=\frac{1}{n^2}\sum_{i=1}^{n}\sum_{j=1}^{n}\frac{1}{4}e^{-2\lambda|j-i|}=\frac{1}{4n^2}\sum_{i=1}^{n}\sum_{j=1}^{n}e^{-2\lambda|j-i|}$$

$$=\frac{1}{4n^2}[n+2\sum_{k=1}^{n-1}(n-k)e^{-2\lambda k}]=\frac{1}{4n^2}[n+2\sum_{k=1}^{n-1}(n-k)e^{-2\lambda k}]$$

$$=\frac{1}{4n^2}[n+2\sum_{k=1}^{n-1}ne^{-2\lambda k}-2\sum_{k=1}^{n-1}ke^{-2\lambda k}]$$

$$=\frac{1}{4n}[1+2\frac{e^{-2\lambda}-e^{-2\lambda n}}{1-e^{-2\lambda}}-\frac{2}{n}\sum_{k=1}^{n-1}ke^{-2\lambda k}]\to 0,\ n\to\infty \tag{33}$$

$$\frac{1}{n}\sum_{k=1}^{n}D(X_k)=\frac{1}{n}\sum_{i=1}^{n}\frac{1}{4}=\frac{1}{4}\not\to 0,\ n\to\infty \tag{34}$$

If we extend the assumption of (1) of *Example* 6 to $P(X(t)=1)=p$ , $P(X(t)=0)=1-p \quad (0<p<1)$ , $\{X(n), n=0,1,2,\cdots\cdots\}$ can be treated as the correlated case of Bernoulli sequence, it is a stationary sequence. Formula (33) shows that the WLLN still holds, while (34) shows that Theorem 7 fails to verify this problem. Hence, we modify Theorem 7 as follows,

**Theorem 12** (Markov Law of large numbers) Let $X_1,\cdots,X_n,\cdots$ be any correlated random variables sequences with $Var(X_i)<+\infty$, $i=1,2,\cdots$. For any $\varepsilon>0$, if

$$(1)\frac{1}{n}\sum_{i=1}^{n}Var(X_i)\to 0 \qquad (n\to+\infty) \tag{35}$$

OR

(2) $\exists\, C>0$, $0<s<1$, satisfying following conditions (a) and (b)

(a). $|\frac{1}{n}\sum_{i=1}^{n}Var(X_i)|<C$ , or $|\frac{1}{n}\sum_{i=1}^{n}Var(X_i)|=O(n^s)$,

$$\text{(b). } |\frac{1}{n}\sum_{1\le i<j\le n}^{n}Cov(X_i,X_j)|<C, \text{ or } |\frac{1}{n}\sum_{1\le i<j\le n}^{n}Cov(X_i,X_j)|=O(n^s) \ . \tag{36}$$

Then,

$$\lim_{n\to\infty}P(|\frac{1}{n}\sum_{i=1}^{n}X_i-\frac{1}{n}\sum_{i=1}^{n}EX_i|<\varepsilon)=1,$$

where $|\frac{1}{n}\sum_{i=1}^{n}Var(X_i)|=O(n^s)$ denotes that $|\frac{1}{n}\sum_{i=1}^{n}Var(X_i)|$ is the same order infinite of $n^s$.

**Proof:**

(1) is the conclusion of Theorem 7.

(2) Since

$$Var(\frac{1}{n}\sum_{i=1}^{n}X_i)=\frac{1}{n}[\frac{1}{n}\sum_{i=1}^{n}Var(X_i)+\frac{2}{n}\sum_{1\le i<j\le n}Cov(X_i,X_j)]$$

When $\exists\, C>0$, $0<s<1$, $|\frac{1}{n}\sum_{i=1}^{n}Var(X_i)|<C$ , or $|\frac{1}{n}\sum_{i=1}^{n}Var(X_i)|=O(n^s)$, and

$|\frac{1}{n}\sum_{1\le i<j\le n}^{n}Cov(X_i,X_j)|<C$, or $|\frac{1}{n}\sum_{1\le i<j\le n}^{n}Cov(X_i,X_j)|=O(n^s)$ hold, we have

$$Var(\frac{1}{n}\sum_{i=1}^{n}X_i)=\frac{1}{n}[\frac{1}{n}\sum_{i=1}^{n}Var(X_i)+\frac{2}{n}\sum_{1\le i<j\le n}Cov(X_i,X_j)]\to 0,\ n\to\infty.$$

Using Chebyshev's inequality,

$$P(|\frac{1}{n}\sum_{i=1}^{n} X_i - \frac{1}{n}\sum_{i=1}^{n} EX_i| \geq \varepsilon) \leq \frac{1}{\varepsilon^2} Var(\frac{1}{n}\sum_{i=1}^{n} X_i) \to 0, \ n \to \infty$$

Which ends the proof of the theorem. □

**Lemma 2.** Let $a_1, a_2, \cdots, a_n, \cdots$ be any real number sequence, if $\lim_{n\to\infty} a_n = a$, then

$$\lim_{n\to\infty} \frac{a_1 + \cdots + a_n}{n} = \lim_{n\to\infty} a_n = a.$$

Consider the condition (1) in Theorem 12, the same distribution sequence absolutely does not satisfy it, according to Lemma 2, we obtain a sufficient condition of condition (1) in Theorem 12 as follows.

**Theorem 13.** Let $X_1, \cdots, X_n, \cdots$ be any correlated random variables sequences, $n = 1, 2, \cdots$. If

(1) $Var(X_n) \to 0 \quad (n \to +\infty)$ (37)

OR

(2) $\exists C > 0$, $0 < s < 1$, satisfying the following conditions (a) and (b)

(a) $|\frac{1}{n}\sum_{i=1}^{n} Var(X_i)| < C$, or $|\frac{1}{n}\sum_{i=1}^{n} Var(X_i)| = O(n^s)$,

(b) $|\frac{1}{n}\sum_{1\le i<j\le n} Cov(X_i, X_j)| < C$, or $|\frac{1}{n}\sum_{1\le i<j\le n} Cov(X_i, X_j)| = O(n^s)$.

Then,

$$\lim_{n\to\infty} P(|\frac{1}{n}\sum_{i=1}^{n} X_i - \frac{1}{n}\sum_{i=1}^{n} EX_i| < \varepsilon) = 1.$$

Note that there is a famous sufficient condition for WLLN ([6],p229; [17],p90-91).

**Lemma 3.** Let $X_1, \cdots, X_n, \cdots$ be any correlated random variables sequence with finite variance $Var(X_n) < C$, if $Cov(X_k, X_l) \to 0$, as $|k - l| \to +\infty$, then,

$$\lim_{n\to\infty} P(|\frac{1}{n}\sum_{i=1}^{n} X_i - \frac{1}{n}\sum_{i=1}^{n} EX_i| < \varepsilon) = 1.$$

Note that the conclusion of Lemma 3 is different from conditions of Theorem 12 and Theorem 13.

Firstly, the conclusion of finite variance $Var(X_n) < C$ can not guarantee $\frac{1}{n}\sum_{i=1}^{n} Var(X_i) \to 0$, $(n \to +\infty)$ or $Var(X_n) \to 0$, $(n \to +\infty)$. Inversely, according to Cauchy-Schwarz inequality of formula (1) , condition $Var(X_n) \to 0$ , $(n \to +\infty)$ , would deduce to $Cov(X_k, X_l) \to 0$ , $|k-l| \to +\infty$. Secondly, the second condition in Theorem 12 and Theorem 13 is weaker than Lemma 3. A counterexample is as follow, where symmetric random variable [1] means that for discrete random variable, the distribution function satisfies $P(X = -x) = P(X = x)$, and for continuous random variable, the density probability function satisfies $f(-x) = f(x)$.

***Example 7.*** Let $X$ be any symmetric random variable with finite variance $0 < Var(X) < C$, and sequence $X_n = (-1)^n X$, $n = 1, 2, 3, \cdots$. Then,

$$Cov(X_k, X_l) = (-1)^{k+l} Var(X) \not\to 0, \ |k-l| \to +\infty.$$

$X_n$ does not satisfy Lemma 3, however, it satisfies the second condition of Theorem 12 and Theorem 13. Hence, it is subject to WLLN.

Then, the discussion of WLLN of famous Chebyshev and Bernoulli Law of large numbers [15,16] could be verified by the second condition of Theorem 12 when we discuss the correlated case.

In fact, Theorem 7,8,9,10,11,12,13 discuss the stable [9, P61-62, 66-67] in the case of correlated random variables (no independence assumption). If any two random variables are uncorrelated (or independent), Theorem 7,8,9,10,11,12,13 still hold, where $s = 0$. Theorem 12 could be applied to the first moment estimation in stationary and non-stationary (while first moment is constant) time series. And, we can also construct the high order moment estimation method by generalize Theorem 12 to discuss $\lim_{n\to\infty} P(|\frac{1}{n}\sum_{i=1}^{n} X_i^r - \frac{1}{n}\sum_{i=1}^{n} EX_i^r| < \varepsilon) = 1$,where $r > 0$, we omit the discussion for the sake of paper length.

In addition, the upper bound in Theorem 1' can be similarly applied to the discussions of Chebyshev inequality and law of large numbers, we omit them here. And, the theoretical discussion will be given in the future work.

Finally, there is a theorem of sufficient and necessary condition of WLLN as follows,

**Lemma 4** Let $X_1,\cdots,X_n,\cdots$ be any correlated random variables sequence, the sufficient and necessary condition of any sequence $\{X_n\}$ obeying WLLN is

$$\lim_{n\to\infty} E\left(\frac{[\sum_{i=1}^{n}(X_n - EX_n)]^2}{n^2 + [\sum_{i=1}^{n}(X_n - EX_n)]^2}\right) = 0 .$$

As this condition is not easy to verify in complex distribution case of $\{X_n\}$, the sufficient condition investigation of this paper still has practical and theoretical significance. At the end of the discussion, we give an example, that could be examined by Lemma 3, Theorem 12,Theorem 13, and Lemma 4.

***Example 8.*** Let $X$ be random variable with $P(X=-1)=P(X=1)=\frac{1}{2}$ , and sequence $X_n = (-1)^n X$ , $n = 1,2,3,\cdots$. Then,

1) $$\frac{[\sum_{i=1}^{n}(X_n - EX_n)]^2}{n^2 + [\sum_{i=1}^{n}(X_n - EX_n)]^2} = \begin{cases} 0, & n = even\ \ number \\ \dfrac{X^2}{n^2 + X^2}, & n = odd\ \ number \end{cases}$$

$$E\left(\frac{[\sum_{i=1}^{n}(X_n - EX_n)]^2}{n^2 + [\sum_{i=1}^{n}(X_n - EX_n)]^2}\right) = \begin{cases} 0, & n = even\ \ number \\ \dfrac{1}{n^2+1}, & n = odd\ \ number \end{cases}$$

$$\lim_{n\to\infty} E\left(\frac{[\sum_{i=1}^{n}(X_n - EX_n)]^2}{n^2 + [\sum_{i=1}^{n}(X_n - EX_n)]^2}\right) = 0 .$$

Hence, according to Lemma 4, $X_n = (-1)^n X$ obeys WLLN.

2) Since $Cov(X_k, X_l) = (-1)^{k+l} Var(X) \not\to 0,\ |k-l| \to +\infty$ , Lemma 3 fails to verify the conclusion of WLLN.

3) Since $|\frac{1}{n}\sum_{1\le i<j\le n}^{n} Cov(X_i, X_j)| = \begin{cases} Var(X), & n = even\, number \\ \dfrac{n-1}{n} Var(X), & n = odd\ \ number \end{cases}$, according to the second condition of Theorem 12 and Theorem 13, it is subject to WLLN.

$X_n$ does not satisfy Lemma 3, however, it satisfies Theorem 12 and Theorem 13. Combining the conditions of Lemma 3, Theorem 12 and Theorem 13 could provide more convenient sufficient conditions in different cases for WLLN.

**Theorem 14.** Let $X_1, \cdots, X_n, \cdots$ be any correlated random variables sequences, $n = 1, 2, \cdots$. If one of the following conditions is satisfied,

(1) $Var(X_n) \to 0 \quad (n \to +\infty)$

(2) $\frac{1}{n}\sum_{i=1}^{n} Var(X_i) \to 0 \quad (n \to +\infty)$

(3) $\exists C > 0$, $0 < s < 1$, satisfying the following conditions (a) and (b)

(a) $|\frac{1}{n}\sum_{i=1}^{n} Var(X_i)| < C$ , or $|\frac{1}{n}\sum_{i=1}^{n} Var(X_i)| = O(n^s)$,

(b) $|\frac{1}{n}\sum_{1 \le i < j \le n}^{n} Cov(X_i, X_j)| < C$, or $|\frac{1}{n}\sum_{1 \le i < j \le n}^{n} Cov(X_i, X_j)| = O(n^s)$ .

(4) $Var(X_n) < C$ , and $Cov(X_k, X_l) \to 0$, $|k - l| \to +\infty$.

Then,

$$\lim_{n \to \infty} P(|\frac{1}{n}\sum_{i=1}^{n} X_i - \frac{1}{n}\sum_{i=1}^{n} EX_i| < \varepsilon) = 1.$$

# 5. Conclusion

The variance inequalities for weighted sum of correlated random variables are established by two methods Cauchy-Schwarz's inequality and positive semidefinite matrix. And the Chebyshev inequality is also extended to sum of correlated random variables. We also discuss the Markov WLLN, Chebyshev WLLN and Bernoulli WLLN in correlated random variables cases. In our discussion, independence and identical distribution are not necessary, and we only discuss the random variables with second moments. The WLLN discuss of correlated random variables could be utilized in the stationary or non-stationary time series, when their first moments are the same constant, the point estimation holds, otherwise the WLLN describes the stable convergence problem. Future work will focus on applying the conclusions to modifying probability inequalities and analyzing statistical data in various applied fields.